\documentclass[leqno]{article}
\usepackage{amssymb}
\usepackage{amsmath}
\usepackage{amsfonts}

\newtheorem{theorem}{Theorem}

\newtheorem{lemma}[theorem]{Lemma}
\newenvironment{proof}[1][Proof]{\noindent\textbf{#1.} }{\ \rule{0.5em}{0.5em}}

\newenvironment{t_enumerate}{
\begin{enumerate}
\setlength{\itemsep}{1pt}
\setlength{\parskip}{0pt}
\setlength{\parsep}{0pt}}{\end{enumerate}
}

\newenvironment{t_itemize}{
\begin{itemize}
\setlength{\itemsep}{1pt}
\setlength{\parskip}{0pt}
\setlength{\parsep}{0pt}}{\end{itemize}
}

\numberwithin{equation}{section}
\numberwithin{theorem}{section}

\begin{document}

\date{}
\title{Higher asymptotics of unitarity in \\``quantization commutes with reduction''}
\author{William D. Kirwin}
\maketitle

\begin{abstract}
Let $M$ be a compact K\"{a}hler manifold equipped with a Hamiltonian action of a compact Lie group $G$. In [Invent. Math. \textbf{67} (1982), no.~3, 515--538], Guillemin and Sternberg showed that there is a geometrically natural isomorphism between the $G$-invariant quantum Hilbert space over $M$ and the quantum Hilbert space over the symplectic quotient $M/\!\!/G$. This map, though, is not in general unitary, even to leading order in $\hslash.$

In [Comm. Math. Phys. \textbf{275} (2007), no.~2, 401--422], Hall and the author showed that when the metaplectic correction is included, one does obtain a map which, while not in general unitary for any fixed $\hslash,$ becomes unitary in the semiclassical limit $\hslash\rightarrow0$ (\textit{cf.} the work of Ma and Zhang in [C. R. Math. Acad. Sci. Paris \textbf{341} (2005), no.~5, 297--302], and [Astérisque No. 318 (2008), viii+154 pp.]). The unitarity of the classical Guillemin--Sternberg map and the metaplectically corrected analogue is measured by certain functions on the symplectic quotient $M/\!\!/G$. In this paper, we give precise expressions for these functions, and compute complete asymptotic expansions for them as $\hslash\rightarrow0$.
\end{abstract}


\section{Introduction.}

Let $M$ be a compact K\"{a}hler manifold with K\"{a}hler form $\omega$. Suppose there exists a Hermitian line bundle $\ell$ with connection with curvature $-i\omega$. For each positive integer $k$, the geometric quantization $\mathcal{H}_{M}^{(k)}$ of $M$ is defined to be the space of holomorphic sections of $\ell^{\otimes k}$. In the context of geometric quantization, $k$ is interpreted as the reciprocal of Planck's constant $\hslash$.

Suppose moreover that $G$ is a compact Lie group, with Lie algebra $\mathfrak{g}$, which acts on $M$ in a Hamiltonian fashion with moment map $\Phi:M\rightarrow\mathfrak{g}^{\ast}$. Under sufficient regularity assumptions, the symplectic quotient $M/\!\!/G$ is again a compact K\"{a}hler manifold; denote the resulting K\"{a}hler form by $\widehat{\omega}$. Assuming that the action of $G$ lifts, the bundle $\ell$ descends to a line bundle $\hat{\ell}\rightarrow M/\!\!/G,$ and the connection descends to one with curvature $-i\widehat{\omega}$.

The space $\mathcal{H}_{M/\!\!/G}^{(k)}$ of holomorphic sections of $\hat{\ell}^{\otimes k}$ is the result of reducing \textit{before} quantizing. On the other hand, one may \textit{first} quantize and \textit{then }reduce, which amounts to considering the space $\left(  \mathcal{H}_{M}^{(k)}\right)^{G}$ of $G$-invariant sections of $\ell^{\otimes k}.$

A classical result of Guillemin and Sternberg \cite{Guillemin-Sternberg} is that there is a natural invertible linear map $A_{k}$ from the ``first quantize then reduce'' space $\left(\mathcal{H}_{M}^{(k)}\right)  ^{G}$ to the ``first
reduce then quantize'' space $\mathcal{H}_{M/\!\!/G}^{(k)}$. From the point of view of quantum mechanics, though, it is not only the vector space structure of the quantization that is important, but also the inner product.

It is known that in general, the Guillemin--Sternberg map $A_{k}$ is \textit{not} unitary, and is not even unitary to leading order as $k\rightarrow\infty$ \cite{Charles06Toeplitz}, \cite{Flude}, \cite{Hall-K}, \cite{Li}, \cite{Ma-MarinescuBOOK}, \cite{Ma-Zhang05}, \cite{Ma-Zhang06}, \cite{Paoletti}. In \cite{Hall-K}, the author and Brian Hall showed that when the so-called metaplectic correction is introduced, one obtains an analogue $B_{k}$ of the Guillemin--Sternberg map which, though still not unitary in general for any fixed $k$, becomes unitary in the semiclassical limit $k\rightarrow\infty$. This can also be obtained as a corollary of the work of Ma and Zhang in \cite{Ma-Zhang05}, \cite{Ma-Zhang06}. It was later shown to be the case in a more general setting by Hui Li \cite{Li}.

The unitarity, or lack thereof, of the map $A_{k}$ (resp. $B_{k}$) can be measured by a certain function $I_{k}$ (resp. $J_{k}$) on the symplectic quotient $M/\!\!/G$, with unitarity achieved at least when $I_{k}$ (resp. $J_{k}$) is identically $1$. One of the main results of \cite{Hall-K} is that
\[
\lim_{k\rightarrow\infty}J_{k}=1,
\]
where the limit is uniform on $M/\!\!/G$. (There is an analogous computation for the limit of the $I_{k}$; see Sections \ref{subsec:main_results} and \ref{subsec:previous} below).

Our main results are explicit expressions for $I_{k}$ and $J_{k}$ (Theorem \ref{thm:densities}) as well as complete asymptotic expansion of the functions $I_{k}$ and $J_{k}$ as $k\rightarrow\infty$ (Theorem \ref{thm:main_asymp}). We state these results precisely in Section \ref{subsec:main_results} below. We should also mention that in \cite{Ma-Zhang05} and \cite{Ma-Zhang06}, Ma and Zhang, as well as Ma and Marinescu in \cite{Ma-MarinescuBOOK}, already show the existence of asymptotic series which are related to that for $I_{k}$ and $J_{k}.$

Inspired by results in equivariant cohomology, for example the Duistermaat--Heckman theorem, which one may prove by computing the leading order asymptotics and then showing that in fact the higher order corrections are zero, one may ask whether there are any cases in which $B_{k}$ is asymptotically unitary to all orders, that is, in which $\lim_{k\rightarrow\infty}J_{k}=1+o(k^{-\infty})$. Although we do not prove it here, our results suggest that such ``exact asymptotics'' are not possible for compact $M$ (see the remark following Lemma \ref{lemma:series}). Our results, as well as those of \cite{Hall-K}, do not seem to depend crucially on the compactness of $M$, and indeed the obstruction to ``exact asymptotics'' disappears when $M$ is noncompact.

\bigskip
In the rest of this section, we describe our main results precisely and then recall from \cite{Hall-K} the precise definition of modified Guillemin--Sternberg-type map $B_{k}$. We finish this section by explaining how the existence of the asymptotic expansions of $I_k$ and $J_k$ can be deduced from the work of Ma and Zhang in \cite{Ma-Zhang05}, \cite{Zhang}, which also yields another method for computing the coefficients. In Section \ref{sec:densities} we build on the results of \cite{Hall-K} to give precise expressions for the densities $I_{k}$ and $J_{k}$ which make the asymptotic computations possible. In Section \ref{sec:proof} we prove our main result, Theorem \ref{thm:main_asymp}, by applying previous results of the author \cite{K-Lap} to the case at hand.

\subsection{Main results.\label{subsec:main_results}}

Let $(M^{2n},\omega,J,B:=\omega(\cdot,J\cdot))$ be a compact K\"{a}hler manifold with symplectic form $\omega,$ complex structure $J$ and metric $B$. Let $\ell\rightarrow M$ be a Hermitian line bundle over $M$ with connection $\nabla$ with  curvature $-i\omega$. Suppose that a compact Lie group $G$ of dimension $d$ acts on $M$ (preserving the K\"{a}hler  structure) in a Hamiltonian fashion with moment map $\Phi:M\rightarrow\mathfrak{g}^{\ast}$, and suppose moreover that the induced infinitesimal action on $\ell$ (given by the quantization of the components of the moment map) exponentiates to a global action of $G$ on $\ell$. We denote the components of the moment map by $\phi_{\xi}:M\rightarrow\mathbb{R}$, for $\xi\in\mathfrak{g}.$

Suppose that $0$ is a value and a regular value of $\Phi$, and moreover that $G$ acts freely on the zero set $\Phi^{-1}(0)$. In this case, the symplectic quotient $M/\!\!/G:=\Phi^{-1}(0)/G$ is a compact smooth manifold which
inherits a K\"{a}hler structure from that of $M$; denote the induced symplectic form on $M/\!\!/G$ by $\widehat{\omega}.$ The line bundle $\ell$ descends to a Hermitian line bundle $\widehat{\ell}\rightarrow M/\!\!/G$, and the connection $\nabla,$ restricted to $G$-invariant sections, induces a connection on $\widehat{\ell}.$ Throughout, $x_{0}\in\Phi^{-1}(0)$ will denote a point in the zero-section, and $[x_{0}]:=G\cdot x_{0}$ will denote the corresponding point in the symplectic quotient.

The infinitesimal action of $G$ on $M$ can be continued to an infinitesimal action of the complexified group\footnote{For each compact Lie group $G$ there exists a unique Lie group $G_{\mathbb{C}}$ such that $G$ is a maximal compact subgroup which sits inside $G_{\mathbb{C}}$ as a totally real submanifold, and such that the Lie algebra of $G_{\mathbb{C}}$ is the complexification of $\mathfrak{g}$. The group $G_{\mathbb{C}}$ is called the complexification of $G_{\mathbb{C}}.$ It is diffeomorphic to $T^{\ast}G$, and the multiplication map $G\times\exp(\sqrt{-1}\mathfrak{g})\rightarrow G_{\mathbb{C}}$ is a diffeomorphism. See \cite{Knapp}, Section VII.1, for details.} $G_{\mathbb{C}}$ by setting $X^{\sqrt{-1}\xi}:=JX^{\xi}.$ This action exponentiates to an action of $G_{\mathbb{C}}$ on $M$. The saturation $G_{\mathbb{C}}\cdot\Phi^{-1}(0)$ of the zero set by the group $G_{\mathbb{C}}$ is called the stable set $M_{s}$. It is an open submanifold of $M$, and the complement is of complex codimension at least one. The (free) action of $G_{\mathbb{C}}$ on $M_{s}$ gives the stable set the structure of a principle $G_{\mathbb{C}}$-bundle $\pi_{\mathbb{C}}:M_{s}\rightarrow M/\!\!/G$. Indeed, the complex structure on the symplectic quotient can be understood via the K\"{a}hler isomorphism
\[
\Phi^{-1}(0)/G=M/\!\!/G=M_{s}/G_{\mathbb{C}}.
\]
Moreover, the action $\Lambda:\exp(\sqrt{-1}\mathfrak{g})\times\Phi^{-1}(0)\rightarrow M_{s}$ gives the stable set the structure of a trivial vector bundle\footnote{This does \textit{not }imply that $M_{s}$ is a trivializable $G_{\mathbb{C}}$-bundle over $M/\!\!/G$; in general the zero set $\Phi^{-1}(0)$ is a non-trivial $G$-bundle over $M/\!\!/G$.} over $\Phi^{-1}(0)$ with fiber $\mathfrak{g}$ (see the original paper by Guillemin and Sternberg, \cite{Guillemin-Sternberg}, or the more recent articles \cite{Hall-K}, \cite{Li}, or \cite{Sjamaar95}, for details).

\bigskip
The geometric quantization $\mathcal{H}_{M}^{(k)}$ of $M$ is the space of holomorphic sections of $\ell^{\otimes k},~k\in\mathbb{N}$. We make $\mathcal{H}_{M}^{(k)}$ into a Hilbert space by equipping it with the inner product
\[
\left\langle s,t\right\rangle :=\left(  k/2\pi\right)  ^{n/2}\int
_{M}(s,t)\frac{\omega^{n}}{n!},
\]
where $(s,t)$ denotes the pointwise Hermitian product in $\ell^{\otimes k}$.

Let $K:\bigwedge^{n}\left(  T^{1,0}M\right)  ^{\ast}$ denote the canonical bundle of $M.$ Suppose that $K$ admits a square root\footnote{A square root of $K$ is a line bundle, denoted by $\sqrt{K}$, such that $\sqrt{K}\otimes \sqrt{K}=K$. Such a line bundle exists if the second Stiefel--Whitney class of $M$ vanishes. If a square root of $K$ exists, then the set of isomorphism classes is parameterized by $H^1(M,\mathbb{Z}_{2})$.}, and denote a choice of square root by $\sqrt{K}.$ Sections of $\sqrt{K}$ are called half-forms, and $\sqrt{K}$ is called a half-form bundle. A section of $K$ is said to be holomorphic if in each local holomorphic coordinate chart, the coefficient of $dz^{1}\wedge\cdots\wedge dz^{n}$ is a holomorphic function. Suppose that the action of $G$ lifts to an action on $\sqrt{K}$ which is compatible with the action on $K$ induced by pushforward.

There is a natural inner product on the space of sections of $\sqrt{K}$: if $\mu,\nu\in\Gamma\left(  \sqrt{K}\right)$, then $\mu^{2}\wedge\bar{\nu}^{2}\in\bigwedge^{2n}T^{\mathbb{C}}M$ is a (complex) volume form. We can trivialize the bundle $\bigwedge^{2n}TM$ by the (global, nowhere vanishing) section $\omega^{n}/n!$, and hence there is a function $(\mu,\nu)$---the pointwise inner product of $\mu$ and $\nu$---such that
\begin{equation}
\mu^{2}\wedge\bar{\nu}^{2}=:(\mu,\nu)^{2}\omega^{n}/n!.
\label{eqn:halfform-pairing}%
\end{equation}

The metaplectic correction, by definition, amounts to considering $\ell^{\otimes k}\otimes\sqrt{K}$; that is, the (half-form) corrected quantization $\widehat{\mathcal{H}}_{M}^{(k)}$ of $M$ is the space of holomorphic sections of $\ell^{\otimes k}\otimes\sqrt{K}.$ The pairing (\ref{eqn:halfform-pairing}) is a special case of the BKS pairing in geometric quantization \cite{Woodhouse}, Section 10.2. It defines a Hermitian form on $\ell^{\otimes k}\otimes\sqrt{K}$ and hence an inner product on $\widehat{\mathcal{H}}_{M}^{(k)}$: for sections $t_{1},t_{2}\in\Gamma(\ell^{\otimes k}\otimes\sqrt{K})$ which are locally represented by $t_{j}(x)=s_{j}(x)\otimes\mu(x)$, we set
\[
(t_{1},t_{2})(x)=(s_{1}(x),s_{2}(x))(\mu,\nu)(x).
\]

\bigskip
Let $\left(  \mathcal{H}_{M}^{(k)}\right)  ^{G}$ denote the space of $G$-invariant holomorphic sections of $\ell^{\otimes k}$, and similarly $\left(\widehat{\mathcal{H}}_{M}^{(k)}\right)^{G}$ the space of $G$-invariant holomorphic sections of $\ell^{\otimes k}\otimes\sqrt{K}.$ The restriction of a $G$-invariant holomorphic section $s\in\left(  \mathcal{H}_{M}^{(k)}\right)^{G}$ to $\Phi^{-1}(0)$ descends to a section of $\widehat{\ell}$ which we denote by $A_{k}s$. In \cite{Guillemin-Sternberg}, Guillemin and Sternberg show that $A_{k}s$ is holomorphic, and moreover that $A_{k}$ is an isomorphism of vector spaces.

As mentioned above, it is known that the map $A_k$ is generically not unitary, and does not even become approximately unitary as $k\rightarrow\infty.$
In Section \ref{thm:modified GS map} below, we will recall from \cite{Hall-K} a similar map $B_{k}:\left(  \widehat{\mathcal{H}}_{M}^{(k)}\right)^{G}\rightarrow\widehat{\mathcal{H}}_{M/\!\!/G}^{(k)}$, for $k$ sufficiently large, relating the quantum spaces in the presence of the metaplectic correction. To define the map $B_{k}$ requires more than just ``restrict and descend'' because two half-forms on $M$ pair to give an $(n,0)$-form on $M$. But two half-forms on the quotient should pair to give an $(n-d,0)$-form, so a mechanism to reduce the degree is needed. The map $B_{k}$ turns out to be essentially a square root of the map ``restrict to $\Phi^{-1}(0)$, contract with the vectors in the directions of the infinitesimal $G$-action, and descend the $G$-invariant result to the quotient''. The map $B_{k}$ is also in general not unitary, but it \textit{does} become approximately unitary as $k\rightarrow\infty$. The asymptotic unitarity of $B_k$ was shown directly in \cite{Hall-K}, but, as we explain in \ref{subsec:previous} below, it also follows from \cite{Ma-Zhang05} and \cite{Zhang}.

To measure the unitarity of the maps $A_{k}$ and $B_{k}$, the author and Brian Hall showed in \cite{Hall-K} that there exist functions $I_{k}\in C^{\infty }(M/\!\!/G)$, and for $k$ sufficiently large functions $J_{k}\in C^{\infty}(M/\!\!/G)$, such that
\begin{align}
\int_{M}\left\vert s\right\vert ^{2}\frac{\omega^{n}}{n!}  &  =\int_{M/\!\!/G}\left\vert A_{k}s\right\vert ^{2}I_{k}\frac{\widehat{\omega}^{n-d}}{(n-d)!},\text{ for every }s\in\left(  \mathcal{H}_{M}^{(k)}\right)
^{G},\text{and}\label{eqn:Ik-def}\\
\label{eqn:Jk-lim}\int_{M}\left\vert r\right\vert ^{2}\frac{\omega^{n}}{n!}  &  =\int_{M/\!\!/G}\left\vert B_{k}r\right\vert ^{2} J_{k}\frac{\widehat{\omega}^{n-d}}{(n-d)!}, \text{ for every }r\in \left(\widehat{\mathcal{H}}_{M}^{(k)}\right)^{G}.
\end{align}

Clearly, $A_{k}$ (resp. $B_{k}$) is unitary if $I(k)$ (resp. $J(k)$) is identically $1$. The main result of \cite{Hall-K} is a direct proof that for each $x_{0}\in\Phi^{-1}(0)$,
\begin{align*}
\lim_{k\rightarrow\infty}I_{k}([x_{0}])  &  =2^{-d/2}vol(G\cdot x_{0}),\text{and}\\
\lim_{k\rightarrow\infty}J_{k}([x_{0}])  &  =1,
\end{align*}
where both limits are uniform. The asymptotic limit of $B_k$ means in particular, that in the presence of the metaplectic correction, quantization commutes unitarily with symplectic reduction in the semiclassical limit. Moreover, in the uncorrected case, if $vol(G\cdot x_{0})$ is not constant, then $A_{k}$ does not converge to (a constant multiple) of a unitary map.

To the best of our knowledge, the first appearance of the asymptotic limit of $I_k$ is a formal computation in the thesis of Flude \cite{Flude}. The limit (\ref{eqn:Ik-def}) was rigorously obtained first by Charles in \cite{Charles06Toeplitz}, when $G$ is abelian, and by Ma and Zhang in \cite{Ma-Zhang05}, when $G$ is nonabelian. Similar computations have also appeared in \cite{Ma-MarinescuBOOK}, and \cite{Paoletti} (we refer the reader to the discussion in \cite{Hall-K} for more details), and most recently in greater generality in \cite{Li}. As we explain below in Section \ref{subsec:previous}, the limit (\ref{eqn:Jk-lim}) of $J_k$ can also be obtained as a corollary of \cite[Thm. 0.10]{Ma-Zhang05} and \cite{Zhang}.

We will find expressions for the densities $I_{k}$ and $J_{k}$ in terms of the geometric data, and compute complete asymptotic expansions for both densities as $k\rightarrow\infty$.

\bigskip
To state our results precisely, fix an $Ad$-invariant inner product on $\mathfrak{g}$ and, with respect to it, and orthonormal basis $\{\xi_{j}\}_{j=1}^{d}$ such that the corresponding Haar measure $dvol_{G}$ on $G$ is normalized to $\int_{G}dvol_{G}=1$. Introduce polar coordinates $\xi =(\rho,\Omega)$ on $\mathfrak{g}$, where
\[
\rho:=\sqrt{\left(\xi^{1}\right)^{2}+\cdots\left(\xi^{d}\right)^{2}}
\]
and $\Omega\in S^{d-1}$ is a point in the unit sphere; in particular, $\xi=\rho\Omega$. The Lie algebra $\mathfrak{g}$ acts on $M$ infinitesimally, and we denote the vector field giving the action of $\xi\in\mathfrak{g}$ by $X^{\xi}\in\Gamma(TM)$.

For a function $f\in C^{1}(M)$, we define its gradient as the image of $df$ under the isomorphism between $T^{\ast}M$ and $TM$ given by the K\"{a}hler metric $B=\omega(\cdot,J\cdot);$ that is, $df=B(\operatorname{grad}f,\cdot)$. The divergence of a vector field $X\in\Gamma(TM)$ is defined by $\operatorname{div}X:=\mathcal{L}_{X}\omega^{n}/\omega^{n}$, where $\mathcal{L}_{X}$ denotes the Lie derivative in the direction of $X$, since for a K\"{a}hler manifold, the Liouville form $\omega^{n}/n!$ corresponds to the Riemannian volume. These are related to the Laplacian by
\begin{equation}
\Delta f=\operatorname{div}\operatorname{grad}f. \label{eqn:lap_def}
\end{equation}

Our first main result is the following.
\begin{theorem}
\label{thm:densities}The densities $I_{k}$ and $J_{k}$ may be expressed as
\begin{equation}
I_{k}([x_{0}])=\left(\frac{k}{2\pi}\right)^{d/2}vol(G\cdot x_{0})^{2}\,\mathfrak{j}_{1}(k,x_{0}),\text{ and} \label{eqn:Ik}
\end{equation}

\begin{equation}
J_{k}([x_{0}])=\left(\frac{k}{\pi}\right)^{d/2}vol(G\cdot x_{0})\mathfrak{j}_{1/2}(k,x_{0}) \label{eqn:Jk}
\end{equation}
where
\begin{equation}
\mathfrak{j}_{a}(k,x_{0}):=\int_{\mathfrak{g}}\, \exp\left\{\int_{0}^{1}-2k\phi_{\xi}(e^{it\xi}x_{0})+a\,\Delta\phi_{\xi} (e^{it\xi}x_{0})dt\right\} d^{d}\xi. \label{eqn:core_density}
\end{equation}
\end{theorem}

Moreover, we will find that $\mathfrak{j}_{a}(k,x_{0})$ (and hence $I_{k}$ and $J_{k})$ admits an entire asymptotic expansion
\[
\mathfrak{j}_{a}(k,x_{0})\sim k^{-d/2}\sum_{j=0}^{\infty}\zeta_{2j}^{(a)}(x_{0})k^{-j}
\]
as $k\rightarrow\infty$, where the coefficients are given explicitly in Theorem \ref{thm:main_asymp} below. The results of \cite{Hall-K} may be interpreted as the statement that $\zeta_{0}^{(1)}=\zeta_{0}^{(1/2)}=\pi^{d/2}vol(G\cdot x_{0})^{-1}$.

Our second main result is a computation of the coefficients $\zeta_{j}.$ The coefficients can be expressed, and computed, more efficiently in terms of certain combinatorial quantities which we will introduce in Section \ref{sec:proof}. We state here our results in a direct form, where the geometric content can be clearly seen. The concise version appears as Theorem \ref{thm:main_concise} in Section \ref{sec:proof}.

\begin{theorem}
\label{thm:main_asymp}For $\mathfrak{j}_{a}(k)$ as defined in
(\ref{eqn:core_density}),
\begin{equation}
\mathfrak{j}_{a}(k,x_{0})=k^{-d/2}\sum_{j=0}^{\infty}\zeta_{2j}^{(a)}
(x_{0})k^{-j}+o(k^{-\infty}),\label{eqn:asymp_series}%
\end{equation}
where the coefficients are given by
\begin{align}
\zeta_{j}^{(a)}(x_{0}) &  =\tfrac{1}{2}\,\Gamma\left(  \tfrac{d+j}{2}\right) \int_{S^{d-1}}\Bigg[\left\vert X^{\Omega}\right\vert ^{-(d+j)} \sum_{m=0}^{j}\frac{a^{j-m}}{(j-m)!}\nonumber\\
&  \qquad\times\sum_{l=1}^{j-m}\sum_{_{\mathcal{P}_{j,l}(\vec{n})}} c(j-m;\vec{n})(\Delta\phi_{\Omega})^{n_{1}}\left(\left(JX^{\Omega}\right) \Delta\phi_{\Omega}\right)^{n_{2}} \cdots \left(\left(JX^{\Omega}\right)^{l}\Delta\phi_{\Omega}\right)^{n_{l}}\label{eqn:thm3}\\
&  \qquad\times\sum_{r=1}^{m}\binom{-\frac{d+j}{2}}{r}\left\vert X^{\Omega}\right\vert ^{-2r}\sum_{\mathcal{Q}_{m,r}(\vec{n})}\frac{2^{r}}{(n_{1}+2)! \cdots(n_{r}+2)!}(JX^{\Omega})^{n_{1}+1}\phi_{\Omega}\cdots(JX^{\Omega})^{n_{r}+1}\phi_{\Omega}\Bigg]_{x_{0}}d\Omega\nonumber
\end{align}
where
\begin{align*}
\binom{\alpha}{r} &  :=\frac{\alpha(\alpha-1)\cdots(\alpha-r+1)}{r!},\\
c(j;n_{1},n_{2},\dots,n_{j}) &  :=\frac{j!}{(1!)^{n_{1}}n_{1}!(2!)^{n_{2}}n_{2}!\cdots(j!)^{n_{j}}n_{j}!},
\end{align*}
the sums in the second and third lines of (\ref{eqn:thm3}) are taken over the
sets%
\begin{align*}
\mathcal{P}_{j,l}(\vec{n}) &  =\{(n_{1},\dots,n_{l})\in\mathbb{Z}_{\geq0}%
^{l}:n_{1}+\cdots+n_{l}=j-l+1\text{ and }n_{1}+2n_{2}+\cdots+l\,n_{l}%
=j\},\text{ and}\\
\mathcal{Q}_{m,r}(\vec{n}) &  =\{(n_{1},\dots,n_{r})\in\mathbb{Z}_{\geq1}%
^{r}:n_{1}+\cdots+n_{r}=m\},
\end{align*}
and empty sums are understood to be $1$.
\end{theorem}

For example, the first two terms are
\begin{align*}
\zeta_{0}^{(a)} &  =\tfrac{1}{2}\,\Gamma\left(  \tfrac{d}{2}\right)
\int_{S^{d-1}}\left\vert X_{x_{0}}^{\Omega}\right\vert ^{-d}d\Omega,\text{
and}\\
\zeta_{2}^{(a)} &  =\tfrac{d}{4}\,\Gamma\left(  \tfrac{d}{2}\right)
\int_{S^{d-1}}\left\vert X_{x_{0}}^{\Omega}\right\vert ^{-(d+2)}%
\Bigg[a(JX_{x_{0}}^{\Omega}\Delta\phi_{\Omega}(x_{0})+a\left(  \Delta
\phi_{\Omega}(x_{0})\right)  ^{2})\\
&  \qquad\qquad-\tfrac{d+2}{2}\frac{\frac{a}{3}\Delta\phi_{\Omega}(x_{0}) (JX_{x_{0}}^{\Omega})^{2}\phi_{\Omega}(x_{0}) + \frac{1}{12}(JX_{x_{0}}^{\Omega})^{3}\phi_{\Omega}(x_{0})}{\left\vert X_{x_{0}}^{\Omega}\right\vert^{2}} +\binom{-\frac{d+2}{2}}{2}\frac{((JX_{x_{0}}^{\Omega})^{2}\phi_{\Omega}(x_{0}))^{2}}{9\left\vert X_{x_{0}}^{\Omega}\right\vert^{4}}\Bigg]d\Omega.
\end{align*}

\noindent\textbf{Remarks.}

\begin{t_enumerate}
\item The function
\[
JX^{\Omega}\phi_{\Omega}(x_{0})=\omega(X^\Omega,JX^\Omega)=\left\vert X_{x_{0}}^{\Omega}\right\vert ^{2}
\]
is strictly positive (since we assume $G$ acts freely on the zero set). We write it as $|X^\Omega|^2$ when we want to emphasize this positivity.

\item By the general theory of \cite{K-Lap}, the first term can be expressed in terms of the determinant $H$ of the Hessian of $2\int_{0}^{1}\phi_{\xi}(e^{it\xi}x_{0})dt$ at $\xi=0$. This determinant was computed in \cite{Hall-K}, Lemma 3.1 and Theorem 4.1, to be $H=2^{d}vol(G\cdot x_{0})^{2}$, from which it follows by \cite{K-Lap}, Proposition 1, and Lemma \ref{lemma:series} that
\[
\zeta_{0}^{(a)}=\tfrac{1}{2}\,\Gamma\left(  \tfrac{d}{2}\right)  \int
_{S^{d-1}}\left\vert X_{x_{0}}^{\Omega}\right\vert ^{-d}d\Omega=\frac
{\Gamma(\frac{d}{2})}{2}\frac{2^{d/2+1}\pi^{d/2}}{\Gamma(\frac{d}{2})\sqrt{H}%
}=\pi^{d/2}vol(G\cdot x_{0})^{-1}.
\]
$\hfill\lozenge$
\end{t_enumerate}

\bigskip We conclude this section by recalling the definition from \cite{Hall-K} of the half-form corrected Guillemin--Sternberg type map $B_{k}$ discussed in the introduction. For a $p$-form $\alpha$, denote the (left) contraction with vector fields $X_{1},\dots,X_{r}$ by $\mathfrak{i}\left(\bigwedge_{j=1}^{r}X_{j}\right)  \alpha:=\alpha(X_{1},X_{2},\dots,X_{r},\cdot,\dots,\cdot).$

\begin{theorem}
\label{thm:modified GS map}\cite{Hall-K} (Theorem 3.1) There exists a linear map
$B:\Gamma(M,\sqrt{K})^{G}\rightarrow\Gamma(M/\!\!/G,\sqrt{\widehat{K}}),$
unique up to an overall sign, with the property that%
\[
\pi_{\mathbb{C}}^{\ast}\left[  (B\nu)^{2}\right]  =\left.  \left[
\mathfrak{i}\left(  \bigwedge\nolimits_{j}X^{\xi_{j}}\right)  (\nu
^{2})\right]  \right\vert _{M_{s}}.
\]
For any open set $U$ in $M/\!\!/G,$ if $\nu$ is holomorphic in a neighborhood $V$ of $\pi_{\mathbb{C}}^{-1}(U),$ then $B\nu$ is holomorphic on $U.$

For each $k,$ there is a linear map $B_{k}:\Gamma(M,\ell^{\otimes k}\otimes\sqrt{K})^{G}\rightarrow\Gamma(M/\!\!/G,\hat{\ell}^{\otimes k}\otimes\sqrt{\widehat{K}}),$ unique up to an overall sign, with the property that
\[
B_{k}(s\otimes\nu)=A_{k}(s)\otimes B(\nu)
\]
for all $s\in\Gamma(\ell^{\otimes k})$ and $\nu\in\Gamma(\sqrt{K}).$ This map takes holomorphic sections of $\left.  \ell^{\otimes k}\otimes\sqrt{K}\right\vert _{V}$ to holomorphic sections of $\left.  \hat{\ell}^{\otimes k}\otimes\sqrt{\widehat{K}}\right\vert _{U}.$
\end{theorem}

\subsection{Previous work.}
\label{subsec:previous}
Let us briefly describe how the asymptotic unitarity of the maps $B_k$, as well as the existence of the asymptotic expansions of $I_k$ and $J_k$, can be deduced from results of Ma and Zhang \cite{Ma-Zhang05},\cite{Zhang}.

In \cite[Prop. 1.2 and Thm 1.1]{Zhang}, it is shown that there is a canonical isomorphism
\begin{equation}
\label{eqn:zhang1}
H^0(M,\ell^{\otimes k}\otimes\sqrt{K})^G\rightarrow H^0(M/\!\!/G,\hat{\ell}^{\otimes k}\otimes(\sqrt{K})_G),\text{ for }k\gg1.
\end{equation}
In \cite{Hall-K}, the author and Brian Hall show that there is a canonical isomorphism $K_G\rightarrow\widehat{K}$, $\alpha\mapsto B(\alpha),$ given by contracting with the vector fields generating the action of $G$. Combining these two isomorphisms shows that $B_k$ is an isomorphism for $k\gg1.$
Moreover, it follows from the isomorphism (\ref{eqn:zhang1}) that the modification of the norms on $(\sqrt{K})_G$ and $\sqrt{\widehat{K}}$ is induced by $B$ and equals $\sqrt{vol(G\cdot x_0)}.$ Hence, the asymptotic unitarity of $J_k$ is a consequence of \cite[Thm. 0.10]{Ma-Zhang05}.

To see that the asymptotic expansion of $I_k$ follows from \cite{Charles06Toeplitz}, when $G$ is abelian, and from \cite{Ma-Zhang05}, when $G$ is nonabelian, let $s\in(\mathcal{H}^{(k)}_M)^G$ and set $v=A_k s.$ Then
\[\int_M |s|^2\frac{\omega^n}{n!}=\langle(A_k\circ A_k^*)^{-1}v,v\rangle_{M/\!\!/G}.\]
By \cite{Charles06Toeplitz} and \cite{Ma-Zhang05}, the operator $k^{-d/2}(A_k\circ A_k^*)$ is a Toeplitz operator, whence $k^{d/2}(A_k\circ A_k^*)^{-1}$ is also a Toeplitz operator. Thus,
\[k^{d/2}(A_k\circ A_k^*)^{-1}=\sum_{j=0}^\infty \Pi_k g_j \Pi_k k^{-j}+O(k^{-\infty}),\]
where $\Pi_k:C^\infty(\Lambda^{0,\bullet}T^*M/\!\!/G\otimes\hat{\ell}^k)\rightarrow \mathcal{H}_{M/\!\!/G}^{(k)}$ is the orthogonal projection, and $g_j$ are smooth functions on $M/\!\!/G$ with $g_0=2^{-d/2}vol(G\cdot x_0).$ By integrating, we obtain
\[\langle(A_k\circ A_k^*)^{-1}v,v\rangle_{M/\!\!/G}=k^{-d/2}\sum_{j=0}^\infty\int_{M/\!\!/G}|A_k s|^2 g_j k^{-j}\frac{\omega^n}{n!}+O(k^{-\infty}),\]
which shows that the desired expansion for $I_k$ exists. A similar computation can be made to conclude the existence of an expansion for $J_k$ (using the operator $\mathbf{T}_p$ from \cite[(7.1)]{Ma-Zhang05}.

\section{The densities $I_{k}$ and $J_{k}$.\label{sec:densities}}

In this section, we will build on the results of \cite{Hall-K} to find the expressions (\ref{eqn:Ik}) and (\ref{eqn:Jk}) for the densities $I_{k}$ and $J_{k}$ (resp.). Much of the groundwork was already done in \cite{Hall-K}, but we need more precise results to get the full asymptotic expansion.

In \cite[Thm. 4.2 and 4.3]{Hall-K}, it was shown that
\begin{equation}
I_{k}([x_{0}])=vol(G\cdot x_{0})(k/2\pi)^{d/2}\int_{\mathfrak{g}}\tau
(\xi,x_{0})\exp\left\{  -2k\int_{0}^{1}\phi_{\xi}(e^{it\xi}x_{0})dt\right\}
d^{d}\xi\label{eqn:Hall-K_I}%
\end{equation}
and that
\begin{equation}
J_{k}([x_{0}])=(k/\pi)^{d/2}\int_{\mathfrak{g}}\tau(\xi,x_{0})\exp\left\{
-\int_{0}^{1}\left(  2k\phi_{\xi}(e^{it\xi}x_{0})+\frac{\mathcal{L}_{JX^{\xi}%
}\omega^{n}}{2\omega^{n}}(e^{it\xi}x_{0})\right)  dt \right\}  d^{d}\xi,
\label{eqn:Hall-K_J}%
\end{equation}
where $\tau$ is the Jacobian of the diffeomorphism $\Lambda:\mathfrak{g}\times\Phi^{-1}(0)\rightarrow M_{s}$ given by $\Lambda(\xi,x_{0}):=e^{i\xi}x_{0}.$ It was also shown that as $k\rightarrow\infty$, the contribution to $I_{k}$ (resp. $J_{k}$) coming from the complement of a ball of finite radius is exponentially small, so it is enough to consider the integrals restricted to the unit ball $B:=\{\xi\in\mathfrak{g}:\left\vert \xi\right\vert \leq1\}$.

The main result of this section is the following computation of the Jacobian $\tau$ in terms of the geometric data.

\bigskip

\begin{theorem}
\label{thm:tau}The Jacobian of the map $\Lambda:(\xi,x_{0})\in\mathfrak{g}%
\times\Phi^{-1}(0)\rightarrow e^{i\xi}x_{0}\in M_{s}$ is given by
\[
\tau(\xi,x_{0})=vol(G\cdot x_{0})\exp\left\{  \int_{0}^{1}\Delta\phi_{\xi
}(e^{it\xi}x_{0})dt\right\}  .
\]

\end{theorem}

Using Theorem \ref{thm:tau} to simplify the densities (\ref{eqn:Hall-K_I}) and (\ref{eqn:Hall-K_J}) above yields Theorem \ref{thm:densities}. The proof of Theorem \ref{thm:tau} depends on the following two technical lemmas. Consider the volume form $\mu$ on $M_{s}$ given by
\[
\mu_{e^{i\xi}x_{0}}:=(e^{-i\xi})^{\ast}\frac{\omega^{n}}{n!}.
\]

\begin{lemma}
\label{lemma:mu_inv} For each $\eta\in\mathfrak{g}$, we have $\mathcal{L}_{JX^{\eta}}\mu=0.$
\end{lemma}

\begin{proof}
Let $X_{1},\cdots,X_{2n}$ be vector fields on $M_{s}$ in a neighborhood of $e^{i\xi}x_{0}.$ Then
\begin{align*}
\left(  \mathcal{L}_{JX^{\eta}}\mu\right)  _{e^{i\xi}x_{0}}(X_{1},\dots
,X_{2n})  &  =\lim_{s\rightarrow0}\frac{1}{s}\left[  (e^{is\eta})^{\ast}\mu_{e^{i\xi}x_{0}}-\mu_{e^{i\xi}x_{0}}\right]  (X_{1},\dots,X_{2n})\\
&  =\lim_{s\rightarrow0}\frac{1}{s}\left[  \mu_{e^{is\eta+i\xi}x_{0}}(e_{\ast}^{is\eta}X_{1},\dots,e_{\ast}^{is\eta}X_{2n})-\mu_{e^{i\xi}x_{0}}(X_{1}%
,\dots,X_{2n})\right] \\
&  =\lim_{s\rightarrow0}\frac{1}{sn!}\Big[\omega_{x_{0}}^{n}(e_{\ast}%
^{-is\eta-i\xi+is\eta}X_{1},\dots,e_{\ast}^{-is\eta-i\xi+is\eta}X_{2n})\\
&  \qquad\qquad\qquad\qquad-\omega_{x_{0}}^{n}(e_{\ast}^{-i\xi}X_{1}%
,\dots,e_{\ast}^{-i\xi}X_{2n})\Big]=0.
\end{align*}

\end{proof}

\begin{lemma}
\label{lemma:mu_pullback} $\left(  \Lambda^{\ast}\mu\right)  _{(\xi,x_{0})}=vol(G\cdot x_{0})d^{d}\xi\wedge dvol_{\Phi^{-1}(0)}.$
\end{lemma}

\begin{proof}
Since both $\Lambda^{\ast}\mu$ and $d^{d}\xi\wedge dvol_{\Phi^{-1}(0)}$ are top dimensional and the latter is nowhere vanishing, there exists a function $h(\xi,x_{0})$ such that $\left(\Lambda^{\ast}\mu\right)_{(\xi,x_{0})}=h(\xi,x_{0})d^{d}\xi\wedge dvol_{\Phi^{-1}(0)}.$

We will show that $h$ is independent of $\xi$, from which we will conclude that $h(x_{0})=vol(G\cdot x_{0})$ by restricting to the zero set, where it is known \cite{Hall-K}, equation (4.17) and Lemma 5.4, that
\[
(\Lambda^{\ast}\mu)_{(0,x_{0})}=(\Lambda^{\ast}\omega^{n}/n!)_{(0,x_{0})}=vol(G\cdot x_{0})d^{d}\xi\wedge dvol_{\Phi^{-1}(0)}.
\]

To show that $h$ is independent of $\xi$, let $\eta\in\mathfrak{g}$. Then
\[
\Lambda_{\ast}(\eta,0)_{e^{i\xi}x_{0}}=\left.  \frac{d}{ds}\right\vert
_{s=0}\Lambda(\xi+s\eta,x_{0})=\left.  \frac{d}{ds}\right\vert _{s=0}%
e^{i(\xi+s\eta)}x_{0}=JX_{e^{i\xi}x_{0}}^{\eta}.
\]
By Lemma \ref{lemma:mu_inv}, for each $\eta\in\mathfrak{g}$,
\[
\left(  \mathcal{L}_{(\eta,0)}\Lambda^{\ast}\mu\right)  _{(\xi,x_{0})}=\left(
\Lambda^{\ast}\mathcal{L}_{JX^{\eta}}\mu\right)  _{(\xi,x_{0})}=0,
\]
which implies $\frac{\partial}{\partial\xi_{j}}h(\xi,x_{0})=0$ for each basis vector $\xi_{j}$; that is, $h=h(x_{0})$.
\end{proof}

\bigskip
\begin{proof}[Proof of Theorem \ref{thm:tau}]
For each $\xi\in\mathfrak{g}$, the $G_{\mathbb{C}}$-action yields a map $e^{-i\xi}:M_{s}\rightarrow M_{s}$ given by $x\mapsto e^{-i\xi}x$. Let $\mu\in\Gamma\left(  \bigwedge^{2n}T^{\ast}M\right)  $ be the volume form defined at each point by
\[
\mu_{e^{i\xi}x_{0}}:=\left(  (e^{-i\xi})^{\ast}\frac{\omega^{n}}{n!}\right)_{e^{i\xi}x_{0}}.
\]

Since $\mu$ is nonvanishing, we can use it to trivialize $\bigwedge^{2n}T^{\ast}M;$ in particular, there is a function $\delta\in C^{\infty}(M)$ such that
\begin{equation}
\frac{\omega^{n}}{n!}=\delta\mu. \label{eqn:delta_def}%
\end{equation}
Differentiating in the direction of $JX^{\xi}$ and dividing by $\omega^{n}/n!$, we obtain (using Lemma \ref{lemma:mu_inv})
\begin{equation}
\frac{\mathcal{L}_{JX^{\xi}}\omega^{n}}{\omega^{n}}=\frac{JX^{\xi}(\delta
)}{\delta}=JX^{\xi}\log\delta. \label{eqn:dlogdelta}%
\end{equation}

Fix a point $\xi=\rho\Omega$ and define a path $\gamma_{\Omega}(t):=e^{it\Omega}x_{0}$ for $t\in\lbrack0,\rho]$. Then $\dot{\gamma}_{\Omega}(t)=JX_{e^{it\Omega}x_{0}}^{\Omega}.$ Hence, $JX^{\Omega}\log\delta=d\log\delta(\dot{\gamma}_{\Omega}).$ Integrating $d\log\delta$ along the path $\gamma_{\Omega}(t)$ therefore yields
\begin{equation}
\delta(e^{i\rho\Omega}x_{0})=\exp\left\{  \int_{0}^{\rho}\frac{\mathcal{L}_{JX^{\Omega}}\omega^{n}}{\omega^{n}}\right\}  \delta(x_{0}).
\label{eqn:delta_comp}
\end{equation}

Now, by definition $\mu_{x_{0}}=\omega_{x_{0}}^{n}/n!$ for $x_{0}\in\Phi^{-1}(0)$ which implies $\delta(x_{0})=1$ for $x_{0}\in\Phi^{-1}(0)$. Combining Lemma \ref{lemma:mu_pullback} with (\ref{eqn:delta_def}) and (\ref{eqn:delta_comp}) yields
\[
\Lambda^{\ast}\omega^{n}/n!=vol(G\cdot x_{0})\exp\left\{  \int_{0}^{\rho}%
\frac{\mathcal{L}_{JX^{\Omega}}\omega^{n}}{\omega^{n}}\right\}  d^{d}\xi\wedge
dvol_{\Phi^{-1}(0)}.
\]

Finally, to complete the proof, observe that by (\ref{eqn:lap_def}), $\mathcal{L}_{JX^{\xi}}\omega^{n}/\omega^{n}=\operatorname{div}\operatorname{grad}\phi_{\Omega}=\Delta\phi_{\Omega}.$
\end{proof}

\section{The expansion.\label{sec:proof}}

In this section, we first introduce some combinatorial objects to simplify the statement of Theorem \ref{thm:main_asymp}. We then recall results of the author \cite{K-Lap} which can, after some computations, be used to arrive at Theorem \ref{thm:main_asymp}. Finally, we will carry out these computations, thus arriving at our proof of Theorem \ref{thm:main_asymp}.

\bigskip
To state Theorem \ref{thm:main_asymp} in a more useful form, we recall here some combinatorial objects related to Bell polynomials; we refer the interested reader to \cite{Comtet}, Chapter 3, for more details. The \textbf{partial Bell polynomials} $\mathcal{B}_{j,l}=\mathcal{B}_{j,l}(x_{1},x_{2},\dots,x_{l})$, combinatorial functions on the set $\{x_{1},\dots,x_{l}\}$ which can be defined in terms of a formal double series expansion, are given explicitly by
\begin{equation}
\mathcal{B}_{j,l}(x_{1},\dots,x_{l})=\sum_{\mathcal{P}_{j}(\vec{n})}%
c(j;\vec{n})\,x_{1}^{n_{1}}x_{2}^{n_{2}}\cdots x_{l}^{n_{l}} \label{eqn:pBell}%
\end{equation}
where
\[
c(j,\vec{n}):=c(j;n_{1},n_{2},\dots,n_{j}):=\frac{j!}{(1!)^{n_{1}}%
n_{1}!(2!)^{n_{2}}n_{2}!\cdots(l!)^{n_{l}}n_{l}!}%
\]
and the sum is taken over the set $\mathcal{P}_{j,l}(\vec{n})$ consisting of all (ordered) $l$-tuples of nonnegative integers $\vec{n}:=(n_{1},n_{2},\dots,n_{l})$ such that $n_{1}+\cdots+n_{l}=j-l+1$ and $n_{1}+2n_{2}+\cdots
l\,n_{l}=j,$ that is,
\begin{equation}
\mathcal{P}_{j,l}(\vec{n})=\left\{  (n_{1},\dots,n_{l})\in\mathbb{Z}_{\geq
0}^{l}:\left.
\begin{array}
[c]{c}%
n_{1}+n_{2}+\cdots+n_{l}=j-l+1\\
n_{1}+2n_{2}+\cdots+l\,n_{l}=j
\end{array}
\right.  \right\}  . \label{eqn:bellset}%
\end{equation}
\newline The partial Bell polynomials are classical combinatorial objects and are known to satisfy many recursion (and other) identities. We will find useful the combinations
\begin{equation}
B_{j}(x_{1},\cdots,x_{j}):=\sum_{l=1}^{j}\mathcal{B}_{j,l}(x_{1},\dots,x_{l}),
\label{eqn:CEBell}%
\end{equation}
which are known as the \textbf{complete exponential Bell polynomials.}\footnote{Note that the sum starts at $1$. This is in contrast to some other definitions in the literature; generally, sums starting at zero are called more simply the complete Bell polynomials.}

Related, though much simpler, are the polynomials $\mathcal{C}_{m,r}=\mathcal{C}_{m,r}(x_{1},x_{2},\cdots)$ defined by
\[
(x_{1}t+x_{2}t^{2}+x_{3}t^{3}+\cdots)^{r}=\sum_{m=r}^{\infty}\mathcal{C}_{m,r}t^{m}.
\]
These polynomials can be computed recursively via the relation
\[
\mathcal{C}_{m,r}={\sum_{j=r-1}^{m-1}x_{m-j}}\mathcal{C}_{j,r-1}
\]
with initial data $\mathcal{C}_{m,1}(x_{1},x_{2},\dots)=x_{m}.$ Alternatively, $\mathcal{C}_{m,r}$ is the sum of all ordered products of $r$ elements of the set $\{x_{1},x_{2},\dots\}$ such that the subscripts add to $m$:
\begin{equation}
\mathcal{C}_{m,r}(x_{1},\dots,x_{m})=\sum_{\mathcal{Q}_{m,r}(\vec{n})} x_{n_{1}}x_{n_{2}}\dots x_{n_{r}}, \label{eqn:cBell}
\end{equation}
where
\[
\mathcal{Q}_{m,r}(\vec{n})=\{(n_{1},\cdots,n_{r})\in\mathbb{Z}_{>0}^{r}
:n_{1}+n_{2}+\cdots+n_{r}=m\}.
\]

To state our main theorem more concisely, let $f,g,h\in C^{\infty}(\mathfrak{g}\times\Phi^{-1}(0))$ be
\begin{align}
f(\rho,\Omega,x_{0})  &  :=2\int_{0}^{\rho}\phi_{\Omega}(e^{it\Omega x_{0}})dt, \nonumber\\
h(\rho,\Omega,x_{0})  &  :=\int_{0}^{\rho}\Delta\phi_{\Omega}(e^{it\Omega}x_{0})dt\text{, and}\label{eqn:fandg}\\
g(\rho,\Omega,x_{0})  &  :=\exp\left\{ah(\rho,\Omega,x_0)\right\} \nonumber
\end{align}
The expression of Theorem \ref{thm:main_asymp} we give below is in terms of the Bell polynomials introduced above and the radial derivatives of $f$, $g$ and $h$ (which are computed below in Lemma \ref{lemma:series}). It will turn out that the leading order behavior of $f$ is quadratic, so we define
\[
f_{j}(\Omega,x_{0}):=\left.  \frac{1}{(j+2)!}\frac{\partial^{j+2}}{\partial\rho^{j+2}}f(\rho,\Omega,x_{0})\right\vert _{\rho=0}.
\]
For $g$ we define the usual Taylor coefficients
\[
g_{j}(\Omega,x_{0}):=\left.  \frac{1}{j!}\frac{\partial^{j}}{\partial\rho^{j}g(\rho,\Omega,x_{0})}\right\vert_{\rho=0}
\]
and similarly for $h$ (we will drop the $\Omega$ and $x_{0}$ dependence to ease notation).

The following is a more concise version of our main Theorem
\ref{thm:main_asymp}.

\begin{theorem}
\label{thm:main_concise}For $\mathfrak{j}_{a}(k)$ as defined in
(\ref{eqn:core_density}),
\[
\mathfrak{j}_{a}(k,x_{0})=k^{-d/2}\sum_{j=0}^{\infty}\zeta_{2j}^{(a)}%
(x_{0})k^{-j}+o(k^{-\infty}),
\]
where the coefficients are given by%
\begin{align}
\zeta_{j}^{(a)}  &  =\tfrac{1}{2}\Gamma\left(  \tfrac{d+j}{2}\right)
\int_{S^{d-1}}\left[  f_{0}^{-(d+j)}(\Omega)\sum_{m=0}^{j}g_{j-m}\sum
_{r=1}^{j}\binom{-\frac{d+j}{2}}{r}\frac{\mathcal{C}_{m,r}(f_{1},\dots,f_{m}%
)}{f_{0}(\Omega)^{r}}\right]  d\Omega\label{eqn:combcoeffs}\\
&  =\tfrac{1}{2}\Gamma\left(  \tfrac{d+j}{2}\right)  \int_{S^{d-1}}\left[
f_{0}^{-(d+j)}\sum_{m=0}^{j}\frac{a^{j-m}}{(j-m)!}B_{j-m}(h_{1},2!h_{2}%
,\cdots,m!h_{m})\sum_{r=1}^{m}\binom{-\frac{d+j}{2}}{r}\frac{\mathcal{C}%
_{m,r}(f_{1},\dots,f_{m})}{f_{0}{}^{r}}\right]  d\Omega,\nonumber
\end{align}
where
\begin{align*}
p!h_{p}  &  =(JX_{x_{0}}^{\Omega})^{p-1}\Delta\phi_{\Omega}(x_{0}),\\
f_{p}  &  =\frac{2}{(p+2)!}(JX_{x_{0}}^{\Omega})^{p+1}\phi_{\Omega}%
(x_{0})=\frac{2}{(p+2)!}(JX_{x_{0}}^{\Omega})^{p}\left\vert X_{x_{0}}^{\Omega
}\right\vert ^{2},
\end{align*}
the polynomial $\mathcal{C}_{m,r}(f_{1},\dots,f_{m})$ (defined in (\ref{eqn:cBell})) is the sum of all ordered products of $r$ elements of the set $\{f_{1},f_{2},\dots,f_{m}\}$ such that the subscripts add to $m$, the polynomial $B_{j-m}$ is the complete Bell polynomial defined by (\ref{eqn:pBell}) and (\ref{eqn:cBell}), and empty sums are understood to be $1$.
\end{theorem}

\noindent\textbf{Remark.}

To obtain naive ``exact asymptotics'', that is, term-by-term cancelation of the tail of the series (\ref{eqn:asymp_series}), we see from (\ref{eqn:combcoeffs}) that it is necessary that $f=2\int_{0}^{\rho}\phi_{\Omega}(e^{i\rho\Omega t}x_{0})dt$ be quadratic in $\rho$. Otherwise the set $\{f_{1},f_{2},\dots,f_{N}\}$ is nontrivial for all $N>1$ which yields nontrivial $\mathcal{C}_{m,r}(f_{1},\dots,f_{m})$ terms at all orders. For compact $M$, it is not possible that $f$ be quadratic, since if it were, then twice differentiating implies $JX^{\Omega}\phi_{\Omega}(e^{i\rho\Omega}x_{0})=\left\vert X_{e^{i\rho\Omega}x_{0}}^{\Omega}\right\vert ^{2}=const.$ But as $\rho\rightarrow\infty$, the path $e^{i\rho\Omega}x_{0}$ approaches a point $x_{\infty}$ which is fixed by $e^{\Omega}$ (see \cite{Lerman}) so that we must rather have $\left\vert X_{e^{i\rho\Omega}x_{0}}^{\Omega}\right\vert ^{2}\rightarrow0$ as $\rho\rightarrow\infty$. $\hfill\lozenge$

\bigskip

Using the linearity of the moment map and $\xi=\rho\Omega$, we have%
\[
\int_{0}^{1}\phi_{\xi}(e^{it\xi}x_{0})dt=\int_{0}^{\rho}\phi_{\Omega
}(e^{it\Omega}x_{0})dt
\]
and
\[
\int_{0}^{1}\Delta\phi_{\xi}(e^{it\xi}x_{0})dt=\int_{0}^{\rho}\Delta
\phi_{\Omega}(e^{it\Omega}x_{0})dt.
\]
Therefore, in terms of the functions $f$ and $g$ defined in (\ref{eqn:fandg}), the density $\mathfrak{j}_{a}$ may be written
\[
\mathfrak{j}_{a}=\int_{\mathfrak{g}}e^{-kf}g\,d^{d}\xi=\int_{\mathfrak{g}%
}e^{-kf}e^{ah}\,d^{d}\xi,
\]
which is, for each fixed $x_{0}\in\Phi^{-1}(0)$, a Laplace type integral. It follows from \cite{Hall-K} that it is enough to consider the integral over the unit ball $B:=\{\xi\in\mathfrak{g}:\left\vert \xi\right\vert \leq1\}$.

For completeness, we quote the result from \cite{K-Lap} which we need to obtain Theorem \ref{thm:main_asymp}. Let $\{\xi^{1},\dots,\xi^{d}\}$ be coordinates on $\mathbb{R}^{d}$. Denote by $S^{d-1}=\{\left\vert \xi\right\vert =1\}\subset\mathbb{R}^{d}$ the unit sphere and introduce polar coordinates $\rho:=\sqrt{(\xi^{1})^{2}+\cdots(\xi^{d})^{2}}$ and $\Omega=\xi/|\xi|\in S^{d-1}.$

Suppose that $R$ is a region in $\mathbb{R}^{d}$ containing $0$ as an interior point, and let $f$ and $g$ be measurable functions on $R$. Suppose $f$ attains its unique minimum of $0$ at $0$. Assume moreover there exists $N>0$ and

\begin{t_itemize}
\item $N+1$ continuous functions $f_{j}(\Omega),~j=0,\dots,N$ with $f_{0}>0$ such that for some $\nu>0$
\begin{equation}
f(\rho,\Omega)=\rho^{\nu}\sum_{j=0}^{N}f_{j}(\Omega)\rho^{j}+o(\rho^{N+\nu})\text{ as }\rho\rightarrow0,\text{ and} \label{eqn:f-series}
\end{equation}

\item $N+1$ functions $g_{j}(\Omega),~j=0,\dots,N$ such that for some $\lambda>0$
\begin{equation}
g(\rho,\Omega)=\rho^{\lambda-d}\sum_{j=0}^{N}g_{j}(\Omega)\rho^{j}
+o(\rho^{N+\lambda-d})\text{ as }\rho\rightarrow0. \label{eqn:g-series}
\end{equation}
\end{t_itemize}

\begin{theorem}
\cite{K-Lap} \label{thm:KLap}With the hypotheses above, there exists an asymptotic expansion
\begin{equation}
\int_{B}e^{-kf}g\,d^{d}x=\sum_{j=0}^{N}\zeta_{j}k^{-(\lambda+j)/\nu
}+o(k^{-(N+\lambda)/\nu}) \label{eqn:Lap}%
\end{equation}
where the coefficients are given by
\[
\zeta_{j}=\tfrac{1}{\nu}\,\Gamma\left(  \tfrac{j+\lambda}{\nu}\right)
\int_{S^{d-1}}\left[  f_{0}^{-(j+\lambda)/\nu}\sum_{m=0}^{j}g_{j-m}\sum
_{r=1}^{m}\binom{-\frac{j+\lambda}{\nu}}{r}f_{m}^{(r)}f_{0}^{-r}\right]
d\Omega,
\]
where $f_{m}^{(r)}=\mathcal{C}_{m,r}(f_{1},\dots,f_{m})$ is the sum\footnote{For example, $\ f_{6}^{(3)}=6f_{1}f_{2}f_{3}+3f_{1}^{2}f_{4}+f_{2}^{3}$,} of all ordered products of $r$ elements of $\{f_{1},f_{2},\dots,f_{m}\}$ such that the subscripts add to $m,$ $\binom{\alpha}{r}:=\alpha(\alpha-1)\cdots(\alpha-r+1)/r!,$ and empty sums are understood to be $1$.
\end{theorem}

To apply Theorem \ref{thm:KLap}, we need asymptotic expansions of $f$ and $g$ near $0$. We will use their Taylor series:

\begin{lemma}
\label{lemma:series}For $f$ and $g$ as defined in (\ref{eqn:fandg}),
\[
f=\rho^{2}\sum_{j=0}^{\infty}\frac{2\rho^{j}}{(j+2)!} \left(JX^{\Omega}\right)^{j+1} \phi_\Omega +o(\rho^{\infty}),~\rho\rightarrow0
\]
and
\begin{align*}
g  &  =\sum_{j=0}^{\infty}\frac{\rho^{j}a^{j}}{j!}B_{j}(\Delta\phi_{\Omega
}(x_{0}),JX_{x_{0}}^{\Omega}\Delta\phi_{\Omega}(x_{0}),\dots,(JX_{x_{0}%
}^{\Omega})^{j-1}\Delta\phi_{\Omega}(x_{0}))+o(\rho^{\infty}),~\rho
\rightarrow0\\
&  =\sum_{j=0}^{\infty}\rho^{j}\left[  \frac{a^{j}}{j!}\sum_{l=0}^{j}%
\sum_{\mathcal{P}_{j,l}(\vec{n})}c(j,\vec{n})\prod_{p=1}^{l}\left(  \left(
JX_{x_{0}}^{\Omega}\right)  ^{p-1}\Delta\phi_{\Omega}(x_{0})\right)  ^{n_{p}%
}\right]  +o(\rho^{\infty}),~\rho\rightarrow0,
\end{align*}
\newline where $B_{j}$ is the complete exponential Bell polynomial defined in (\ref{eqn:CEBell}) and $\mathcal{P}_{j,l}(\vec{n})$ is defined in
(\ref{eqn:bellset}).
\end{lemma}

\begin{proof}
The fundamental theorem of calculus yields
\[
\partial_{\rho}^{n}f(e^{i\rho\Omega}x_{0})=2\partial_{\rho}^{n-1}%
\phi_{\Omega}(e^{i\rho\Omega}x_{0})=2(JX^{\Omega})^{n-1}\phi_{\Omega
}(e^{i\rho\Omega}x_{0}).
\]
Moreover, $f(x_{0})=\partial_{\rho}f(x_{0})=0$, so that
\[
f=\sum_{j=2}^{\infty}\frac{2\rho^{j}}{j!}\left[(JX^{\Omega})^{j-1}\phi_\Omega\right]_{x_{0}}
\]
as desired.

To compute the Taylor series for $g$, we first recall that the Taylor series for $\exp(h(\rho))$ near $\rho=0$ can be expressed using Fa\`{a} di Bruno's formula as \cite{Comtet}, Section 3.4,
\begin{align}
&  \sum_{j=0}^{\infty}\frac{\rho^{j}}{j!}\exp(h(0))B_{j}(h'(0),\dots,h^{(j)}(0))\label{eqn:exp_taylor}\\
&  =\sum_{j=0}^{\infty}\frac{\rho^{j}}{j!} \left[\exp(h(0))\sum_{l=1}^{j}\sum_{\mathcal{P}_{j,l}(\vec{n})}c(j;\vec{n})\prod_{p=1}^{l} \left(h^{(p)}(0)\right)^{n_{p}}\right],\nonumber
\end{align}
where $B_{j}$ is the complete exponential Bell polynomial (\ref{eqn:CEBell}), and $\mathcal{P}_{j,l}(\vec{n})$ is defined in (\ref{eqn:bellset}). Taking $h(\rho)=a\int_{0}^{\rho}\Delta\phi_{\Omega}(e^{i\Omega t}x_{0})dt$ in
(\ref{eqn:exp_taylor}) and using $\partial_{\rho}g(e^{i\rho\Omega}x_{0})=JX^{\Omega}g(e^{i\rho\Omega}x_{0}),$ $\partial_{\rho}^{(l)}h=(JX^\Omega)^{l-1}\Delta\phi_\Omega,$ and $g(x_{0})=\exp(h(0))=1$
completes the proof.
\end{proof}

\bigskip

We are now ready to prove our main Theorems \ref{thm:main_asymp} and
\ref{thm:main_concise}.

\begin{proof}
[Proof of Theorems \ref{thm:main_asymp} and \ref{thm:main_concise}]Take
$\nu=2$ and $\lambda=d$. From Lemma \ref{lemma:series}, we see that
\begin{equation}
f_{j}(\Omega,x_{0})=\frac{2}{(j+2)!}(JX_{x_{0}}^{\Omega})^{j}\left\vert
X_{x_{0}}^{\Omega}\right\vert ^{2} \label{eqn:fcoeff}%
\end{equation}
and
\begin{align}
g_{j}(\Omega,x_{0})  &  =\frac{a^{j}}{j!}B_{j}(\Delta\phi_{\Omega}%
(x_{0}),JX_{x_{0}}^{\Omega}\Delta\phi_{\Omega}(x_{0}),\dots,(JX_{x_{0}%
}^{\Omega})^{j}\Delta\phi_{\Omega}(x_{0}))\label{eqn:gcoeff}\\
&  =\frac{a^{j}}{j!}\sum_{l=0}^{j}\sum_{\mathcal{P}_{j,l}(\vec{n})}c(j,\vec
{n})\prod_{p=1}^{l}\left(\left(JX_{x_{0}}^{\Omega}\right)  ^{p-1}%
\Delta\phi_{\Omega}(x_{0})\right)^{n_{p}}.\nonumber
\end{align}
Plugging these into Theorem \ref{thm:KLap} applied to $\mathfrak{j}_{a}(k,x_{0})\sim\int_{B}e^{-kf}g\,d^{d}\xi$ yields
\[
\mathfrak{j}_{a}(k,x_{0})=\sum_{j=0}^{\infty}\zeta_{j}k^{-(j+d)/2}+o(k^{-\infty})
\]
where
\begin{equation}
\zeta_{j}=\tfrac{1}{2}\,\Gamma\left(  \tfrac{j+d}{2}\right)  \int_{S^{d-1}
}\left[  \left\vert X_{0}^{\Omega}\right\vert ^{-(d+j)}\sum_{m=0}^{j}
\frac{a^{j-m}B_{j-m}}{(j-m)!}\sum_{r=1}^{m}\binom{-\frac{j+d}{2}}{r}\left\vert
X_{x_{0}}^{\Omega}\right\vert ^{-2r}f_{m}^{(r)}\right]  d\Omega
\label{eqn:asymp1}%
\end{equation}
in which
\begin{align}
f_{m}^{(r)}  &  =\mathcal{C}_{m,r}(\tfrac{2}{3!}\,(JX_{x_{0}}^{\Omega})^2\phi_\Omega, \tfrac{2}{4!}\,(JX_{x_{0}}^{\Omega})^{3}\phi_\Omega,\dots, \tfrac{2}{(m+2)!}\,(JX_{x_{0}}^{\Omega})^{m+1}\phi_\Omega)\label{eqn:fsub}\\
&  =\sum_{\mathcal{Q}_{m,r}(\vec{n})}\frac{2^{r}}{(n_{1}+2)! \cdots(n_{r}+2)!}(JX_{x_{0}}^{\Omega})^{n_{1}+1}\phi_\Omega\cdots(JX^{\Omega})^{n_{r}+1}\phi_\Omega\nonumber
\end{align}
is the sum of all ordered products of $r$ terms of the set $\{f_{1}(\Omega,x_{0}),f_{2}(\Omega,x_{0}),\dots\}$ whose subscripts add to $m$ and
\begin{align}
B_{j-m}  &  =B_{j-m}(\Delta\phi_{\Omega}(x_{0}),JX_{x_{0}}^{\Omega}\Delta
\phi_{\Omega}(x_{0}),\dots,(JX_{x_{0}}^{\Omega})^{j-m-1}\Delta\phi_{\Omega
}(x_{0}))\label{eqn:bsub}\\
&  =\sum_{l=0}^{j-m}\sum_{\mathcal{P}_{j-m,l}(\vec{n})}c(j-m,\vec{n}%
)\prod_{p=1}^{l}\left(  \left(  JX_{x_{0}}^{\Omega}\right)  ^{p-1}\Delta
\phi_{\Omega}(x_{0})\right)  ^{n_{p}}.\nonumber
\end{align}

Aside from simply substituting $f_{m}^{(r)}$ and $B_{j-m}$, we can make one significant simplification:

\begin{lemma}
\label{lemma:jodd}$\zeta_{j}=0$ for $j$ odd.
\end{lemma}

\begin{proof}
The linearity of the infinitesimal action of $G_{\mathbb{C}}$ on $M$ implies $X^{-\Omega}=-X^{\Omega}$ and $JX^{-\Omega}=-JX^{\Omega}.$ Related is the linearity of the moment map in the component index: $\phi_{-\Omega}=-\phi_{\Omega}.$ These facts together imply that $f_{j}(-\Omega,x_{0})=(-1)^{j}f_{j}(\Omega,x_{0}).$ Since the sum of the subscripts of the terms appearing in $f_{m}^{(r)}$ is $m$, we have $f_{m}^{(r)}(-\Omega ,x_{0})=(-1)^{m}f_{m}^{(r)}(\Omega,x_{0})$. Finally, we conclude that
\[
g_{j-m}(-\Omega,x_{0})f_{m}^{(r)}(-\Omega,x_{0})=(-1)^{j}g_{j-m}(\Omega
,x_{0})f_{m}^{(r)}(\Omega,x_{0})
\]
which implies that for $j$ odd, the integrand appearing in $\zeta_{j}$ is antisymmetric with respect to $\Omega\mapsto-\Omega$ so that the integral is $0$ for $j$ odd.
\end{proof}

\bigskip

Making the substitutions of (\ref{eqn:fsub}) and (\ref{eqn:bsub}) into (\ref{eqn:asymp1}) and replacing $j$ by $2j$ (Lemma \ref{lemma:jodd}) yields Theorems \ref{thm:main_asymp} and \ref{thm:main_concise}.
\end{proof}

\subsection*{Acknowledgments}

The author would like to thank the referee for remarks with clarified the exposition, in particular the relation of this work with that of Ma and Zhang in \cite{Ma-Zhang05} and \cite{Ma-Zhang06}. The author is also grateful to the Max Planck Institute for Mathematics in the Sciences for their hospitality during the preparation of this article.

\bigskip

{\small
\providecommand{\bysame}{\leavevmode\hbox to3em{\hrulefill}\thinspace}
\providecommand{\MR}{\relax\ifhmode\unskip\space\fi MR }
\providecommand{\MRhref}[2]{%
  \href{http://www.ams.org/mathscinet-getitem?mr=#1}{#2}
}
\providecommand{\href}[2]{#2}

}

\vfill
\noindent\textsc{Center for Mathematical Analysis, Geometry and Dynamical Systems},
Instituto Superior T\'ecnico, Av. Rovisco Pais, 1049-001, Lisbon, Portugal.\\
Email: will.kirwin@gmail.com

\end{document}